\documentclass[a4paper, leqno, 11pt]{amsart}
\usepackage{amssymb}
\usepackage{graphicx}
\begin{document}
\newcommand*{\threesim}{\mathrel{\vcenter{\offinterlineskip\hbox{$\sim$}\vskip-.35ex\hbox{$\sim$}\vskip-.35ex\hbox{$\sim$}}}}
\newtheorem{defi}{Definition}[section]
\newtheorem{thm}[defi]{Theorem}
\newtheorem{prop}[defi]{Proposition}
\newtheorem*{exa}{Example}
\newtheorem{lem}[defi]{Lemma}
\newtheorem*{thmA}{Theorem A}
\newtheorem*{thmB}{Theorem B}
\newtheorem*{cor}{Corollary}
\theoremstyle{definition}
\newtheorem*{rmk}{Remark}
\theoremstyle{remark}
\newtheorem*{pf}{Proof}
\numberwithin{equation}{section}
\title{On the number of fixed points of sofic flip systems}
\author{Young-One Kim \and Sieye Ryu}
\begin{abstract}
In the case when $X$ is a sofic shift and $\varphi : X \to X$ is a homeomorphism such that $\varphi^2 = \text{id}_X$ and $\varphi \sigma_X = \sigma_X^{-1} \varphi$, 
the number of points in $X$ that are fixed by $\sigma_X^m$ and $\sigma_X^n \varphi$, $m=1,2,\dots$, $n\in\Bbb Z$, is expressed in terms of a finite number of square matrices: 
The matrices are obtained from Krieger's joint state chain of a sofic shift which is conjugate to $X$.
\end{abstract}
\maketitle
\section{Introduction}
Let $(X, T)$ be an invertible topological dynamical system. 
A homeomorphism $F: X \rightarrow X$ is called a \textit{flip map} (or simply a \textit{flip}) for $(X, T)$ if 
\begin{equation}
FT=T^{-1}F \quad \text{and} \quad F^2=\text{id}_X.
\end{equation}
In this case, we call the triple $(X, T, F)$ a flip system. 
Two flip systems $(X, T, F)$ and $(X', T', F')$ are said to be \textit{conjugate} if there is a homeomorphism 
$\theta : X \rightarrow X'$ such that
$$\theta \circ T = T' \circ \theta \quad \text{and} \quad \theta \circ F = F' \circ \theta.$$
We call the homeomorphism $\theta$ a \textit{conjugacy} from $(X, T, F)$ to $(X', T', F').$

If $m$ is a positive integer, the number of periodic points in $(X, T)$ of period $m$ is denoted by $p_m(T)$:
$$p_m(T)=|\{x\in X : T^m(x)=x\}|.$$
Similarly, if $F$ is a flip for $(X, T)$, $m$ is a positive integer and $n$ is an integer, 
then $p_{m, n}(T, F)$ will denote the number of points in $X$ that are fixed by $T^m$ and $T^nF$:
$$p_{m, n}(T, F)=|\{ x\in X : T^m(x)=T^n F(x)=x \}|.$$
It is obvious that if the flip systems $(X, T, F)$ and $(X', T', F')$ are conjugate, 
then $p_m(T)=p_m(T')$ and $p_{m, n}(T, F)=p_{m,n}(T', F')$ for all $m$ and $n$. 
From (1.1), it follows that
$$p_{m,n}(T, F)=p_{m, n+m}(T, F)=p_{m, n+2}(T, F).$$
Hence $p_{m, n}(T, F)=p_{m,0}(T, F)$ if $m$ and $n$ are even; 
$p_{m, n}(T, F)=p_{m, 1}(T, F)$ if $m$ is even and $n$ is odd;
and $p_{m, n}(T, F)=p_{m,0}(T, F)$ for all $n$ if $m$ is odd.

It is well known that if $X$ is a shift of finite type, then there is a square matrix $A$ 
whose entries are non-negative integers such that
\begin{equation}
p_m(\sigma_X)=\text{tr}(A^m)\qquad(m=1, 2, \cdots).
\end{equation}
It is also well known that if $X$ is a sofic shift, then there are square matrices $A_1, A_2, \cdots, A_r$
whose entries are integers such that
$$p_m(\sigma_X)=\sum_{k=1}^{r} (-1)^{k+1} \text{tr}(A_k^m)\qquad(m=1, 2, \cdots).$$

Suppose that the sequence $\{ (p_m(T))^{1/m} \}$ is bounded. 
Then the Artin-Mazur zeta function $\zeta_T$ of $(X, T)$ is defined by 
$$\zeta_T(t)=\exp\left(\sum_{m=1}^\infty \frac{p_m(T)}{m}t^m\right);$$
and the results mentioned above imply that if $X$ is a shift of finite type, 
then $\zeta_{\sigma_X}$ is the reciprocal of a polynimial, and if $X$ is a sofic shift, then $\zeta_{\sigma_X}$ a rational function.

In [KLP], Y.-O. Kim, J. Lee and K.K.Park introduced a zeta function $\zeta_{T, F}$ of a flip system which is given by
$$\zeta_{T, F}(t)=\zeta_T(t^2)^{1/2}\exp(G_{T, F}(t)),$$
where
$$G_{T, F}(t)=\sum_{m=1}^\infty \left(p_{2m-1, 0}(T, F)t^{2m-1}+\frac{p_{2m,0}(T, F)+p_{2m,1}(T, F)}{2}t^{2m}\right).$$
The function $G_{T, F}$ will be called the \textit{generating function} of $(X, T, F).$

In the case when $X$ is a shift of finite type and $\varphi$ is a flip for $(X, \sigma_X)$, 
a formula for $p_{m, n}(\sigma_X, \varphi)$ which is similar to (1.2) is known.
In order to present it, we need some notation: If M is a matrix $\mathcal{S}[M]$ will denote the sum of the entries of M, that is, 
$$\mathcal{S}[M]=\sum_{I, J} \, M_{IJ}$$ 
and $M^{\Delta}$ will denote the matrix whose diagonal entries are identical with those of $M$ but the other entries are equal to 0, that is,
$$(M^{\Delta})_{IJ}=\begin{cases}M_{IJ} \qquad \text{if} \,\,\, I=J, \\ 0 \qquad \quad \,\,\, \text{otherwise}. \end{cases}$$
For instance, $\text{tr}(M) = \mathcal{S}[M^{\Delta}]$, whenever $M$ is a square matrix.

\begin{thmA}
If $X$ is a shift of finite type and $\varphi$ is a flip for $(X, \sigma_X)$, then there are zero-one square matrices $A$ and $J$
such that
\begin{eqnarray*}
p_m(\sigma_X)&=&\mathrm{tr}(A^m),\\
p_{2m,0}(\sigma_X, \varphi) &=& \mathcal{S}[J^{\Delta} A^m J^{\Delta}],\\
p_{2m,1}(\sigma_X, \varphi) &=& \mathcal{S}[(JA)^{\Delta} A^{m-1}(AJ)^{\Delta}] \quad \text{and} \\
p_{2m-1,0}(\sigma_X, \varphi) &=& \mathcal{S}[J^{\Delta}A^{m-1}(AJ)^{\Delta}]\qquad(m=1, 2, \cdots).
\end{eqnarray*}
\end{thmA} 
\begin{rmk}
As it is proved in [KLP], Theorem A implies that if $X$ is a shift of finite type and $\varphi$ is a flip for $(X, \sigma_X)$,
then the generating function of $(X, \sigma_X, \varphi)$ is a rational function.
\end{rmk}
\noindent
In this paper, we extend Theorem A to the case when $X$ is a sofic shift.
\begin{thmB}
If $X$ is a sofic shift and $\varphi$ is a flip for $(X, \sigma_X)$, 
then there are square matrices $A_k, B_k$ and $J_k, \; k=1, 2, \cdots, r,$
whose entries are integers such that
\begin{eqnarray*}
p_m(\sigma_X)&=& \sum_{k=1}^r \, (-1)^{k+1}\, \mathrm{tr}(A^m),\\
p_{2m,0}(\sigma_X, \varphi) &=& \sum_{k=1}^r \, (-1)^{k+1} \, \mathcal{S}\big[J_k^{\Delta} B_k^m J_k^{\Delta}\big],\\
p_{2m,1}(\sigma_X, \varphi) &=& \sum_{k=1}^r \, (-1)^{k+1} \, \mathcal{S}\big[(J_k A_k)^{\Delta} B_k^{m-1}(A_k J_k)^{\Delta}\big] \quad \text{and} \\
p_{2m-1,0}(\sigma_X, \varphi) &=& \sum_{k=1}^r \, (-1)^{k+1} \, \mathcal{S}\big[J_k^{\Delta} B_k^{m-1}(A_k J_k)^{\Delta}\big]\qquad(m=1, 2, \cdots).
\end{eqnarray*}
\end{thmB} 
\noindent
The same method as in the proof of Theorem 3.2 in [KLP] gives the following.
\begin{cor}
If $X$ is a sofic shift and $\varphi$ is a flip for $(X, \sigma_X)$, then the generating function of $(X, \sigma_X, \varphi)$ is a rational function. 
\end{cor}
We prove Theorem B in Section 2. In our proof of the theorem, Krieger's joint state chain plays a crucial role 
and it is described in Section 3. We conclude the paper with an example and some remarks(Section 4).
\section{Proof of Theorem B}
We start with some preliminaries. Suppose $\mathcal{A}$ is a finite alphabet and $\tau:\mathcal{A} \rightarrow \mathcal{A}$ satisfies 
$\tau ^2 =\text{id}_{\mathcal{A}}$. 
Then the map $\varphi_{\tau}:\mathcal{A}^{\mathbb{Z}} \rightarrow \mathcal{A}^{\mathbb{Z}}$ defined by
$$\varphi_{\tau}(x)_i=\tau(x_{-i})$$
is a flip for $(\mathcal{A}^{\mathbb{Z}}, \sigma)$.
If $X$ is a shift space over $\mathcal{A}$ and $\varphi_{\tau}(X)=X$, the restriction of $\varphi_{\tau}$ to $X$ is a flip for $(X, \sigma_X)$.
Conversely, if $X$ is a shift space over $\mathcal{A}$, $\varphi$ is a flip for $(X, \sigma_X)$ and $\varphi$ satisfies 
$$x, x' \in X \;\; \text{and} \;\; x_0=x'_0 \quad \Rightarrow \quad \varphi(x)_0=\varphi(x')_0,$$
then there is a map $\tau:\mathcal{A} \rightarrow \mathcal{A}$ such that 
$\tau^2=\text{id}_{\mathcal{A}}, \; \varphi_{\tau}(X)=X$ and $(\varphi_{\tau})|_X=\varphi$.
In this case, we say that $\varphi$ is a \textit{one-block flip} and we call the restriction of $\tau$ to $\mathcal{B}_1(X)$ the \textit{symbol map} of $\varphi$.

If $A$ is a zero-one $\mathcal{A} \times \mathcal{A}$ matrix, 
$\textsf{X}_A$ will denote the topological Markov chain determined by $A$:
$$\textsf{X}_A=\{ x\in \mathcal{A}^{\mathbb{Z}} \,:\, \forall i\in \mathbb{Z} \;\; A(x_i, x_{i+1})=1 \}.$$
We denote by $\sigma_A$ the restriction of the shift map to $\textsf{X}_A$. 
Let $A$ be such a matrix.
Then we have $\varphi_{\tau}(\textsf{X}_A)=\textsf{X}_A$ whenever
\begin{equation}
A(a,b)=A(\tau(b), \tau(a))\qquad(a, b \in \mathcal{A})
\end{equation}
holds. If we define the $\mathcal{A} \times \mathcal{A}$ matrix $J$ by
\begin{equation}
J(a, b)=\begin{cases} 1 \qquad \text{if} \; \tau(a)=b, \\ 0 \qquad \text{otherwise}, \end{cases}
\end{equation}
then $J$ is a symmetric permutation matrix, $J^2=I$ and the condition (2.1) is equivalent to the equation
\begin{equation}
JA=A^{\textsf{T}}J.
\end{equation}
Conversely, if $J$ is a zero-one $\mathcal{A} \times \mathcal{A}$ matrix such that $J^2=I$, then it is a symmetric permutation matrix
and hence there is a unique map $\tau:\mathcal{A} \rightarrow \mathcal{A}$ 
such that (2.2) holds and consequently $\tau^2=\text{id}_{\mathcal{A}}$.
In this case, the flip $\varphi_{\tau}$ will also be denoted by $\varphi_J$. 
If (2.3) holds, the flip $\varphi_J$ to $\textsf{X}_{A}$ is a flip for $(\textsf{X}_A, \sigma_{A})$.
We denote the restriction by $\varphi_{J,A}$.

The following proposition will be proved in the next section.
\begin{prop}
Suppose $X$ is a sofic shift and $\varphi$ is a flip for $(X, \sigma_X)$. 
Then there are a finite set $\mathcal{A}$, a map $\mathcal{L}:\mathcal{A} \rightarrow \mathcal{B}_1(X)$ 
and zero-one $\mathcal{A} \times \mathcal{A}$ matrices $A$ and $J$ having the following properties:
\newline
\textsc{(1)} $\mathcal{L}_{\infty} : \textsf{X}_A \rightarrow X$ is a factoring.
\newline
\textsc{(2)} $\mathcal{L}_{\infty}$ has no graph diamonds.
\newline
\textsc{(3)} $JA=A^{\textsf{T}}J$ and $J^2=I$.
\newline
\textsc{(4)} $\mathcal{L}_{\infty} \circ \varphi_{J,A} = \varphi \circ \mathcal{L}_{\infty}$.
\newline
\textsc{(5)} If $\delta \in \{0, 1 \}, \, x\in X$ and $\sigma_X^{\delta} \, \varphi(x) =x$, then there is a $y \in \textsf{X}_A$
such that $\mathcal{L}_{\infty}(y)=x$ and $\sigma_A ^{\delta} \, \varphi_{J,A}(y)=y$. 
\end{prop}

In the rest of this section, we suppose that $X$ is a sofic shift, $\varphi$ is a flip for $(X, \sigma_X)$ 
and the system $(\mathcal{A}, \mathcal{L}, A, J)$ is as in the above proposition. We denote by $\tau$ the map of $\mathcal{A}$ onto itself
such that $J(a, \tau(a))=1$ for all $a\in \mathcal{A}$.

If $\pi$ is a permutation of a finite set, its sign will be denoted by $\text{sgn}(\pi)$.
When $<$ is a linear order of $\mathcal{A}, S \subset \mathcal{A}$ and $f:S \rightarrow \mathcal{A}$ is one-to-one,
we set
$$\text{sgn}_<(f)=(-1)^{N(<, \,f)},$$
where
$$N(<,\, f)=|\{ (a,b) \in S \times S : a<b \; \text{and} \; f(b)<f(a) \}|.$$
If $g:f(S) \rightarrow \mathcal{A}$ is one-to-one again, we have
$$\text{sgn}_<(g\circ f)=\text{sgn}_<(g) \, \text{sgn}_<(f).$$
We also have $\text{sgn}_<(f)=\text{sgn}(f)$ wherever $f$ is a permutation of a subset of $\mathcal{A}$.

From here on, we fix a linear order $<$ of $\mathcal{A}$ and write
$$\text{sgn}(f)=\text{sgn}_<(f)$$
for arbitrary one-to-one function $f$ from a subset of $\mathcal{A}$ into $\mathcal{A}$.

Let $k$ be a positive integer $\leq |\mathcal{A}|$. We set
$$\mathcal{A}_k=\{ S \subset \mathcal{A} : |S|=k \; \text{and} \; |\mathcal{L}(S)|=1\} ; $$
and for $S_1, S_2 \in \mathcal{A}_k$ we denote by $F(S_1, S_2)$ the set of one-to-one functions
$f:S_1 \rightarrow S_2$ such that $A(a, f(a))=1$ for all $a\in S_1$.
Define the $\mathcal{A}_k \times \mathcal{A}_k$ matrices $A_k, B_k$ and $J_k$ by
\begin{eqnarray*}
A_k(S_1, S_2) &=& \sum_{f\in F(S_1, S_2)} \text{sgn}(f),\\
B_k(S_1, S_2) &=& \begin{cases} 1 \qquad \text{if} \; F(S_1, S_2)\neq\varnothing, \\
0 \qquad \text{otherwise},\end{cases}\\
\text{and} \hspace{2.5cm}&&\\
J_k(S_1, S_2)&=&\begin{cases} \text{sgn}(\tau|_{S_1}) \quad \text{if} \; \tau(S_1)=S_2,\\
0 \qquad \qquad \; \: \text{otherwise.}\end{cases}
\end{eqnarray*} 

It is well known that
\begin{equation}
p_m(\sigma_X)=\sum_{k=1}^{|\mathcal{A}|} \, (-1)^{k+1} \, \text{tr}(A_k^m) \qquad (m=1,2,\cdots).
\end{equation}
Hence the first formula in Theorem B holds. A proof of (2.4) is found in Section 6.4 of [LM](See also [B] and [M].)
and our proof of the other formulas in Theorem B is a `modification' of it. In order to explain the modification, 
we give an outline of the proof of (2.4).

Let $m$ be a positive integer. We put
$$P(m)=\{ x\in X : \sigma_X(x)=x \},$$
so that $p_m(\sigma_X)=|P(m)|$.
Suppose $x \in P(m)$.
Since $\mathcal{L}_{\infty} : \textsf{X}_A \rightarrow X$ is a factoring with no graph diamonds and since
$\sigma_X^m(x)=x, \, \mathcal{L}_{\infty}^{-1}(x)$ is a non-empty finite subset of $\textsf{X}_A$
and the resriction of $\sigma_A^m$ to $\mathcal{L}_{\infty}^{-1}(x)$ is a permutation of $\mathcal{L}_{\infty}^{-1}(x)$.
Let $\mathcal{C}(m;\,x)$ denote the set of all subsets of $\mathcal{L}_{\infty}^{-1}(x)$ that are fixed by $\sigma_A^m$:
$$\mathcal{C}(m; \, x) = \{ E \subset \mathcal{L}_{\infty}^{-1}(x) : \sigma_A^m(E)=E \}.$$

The following lemma is proved in section 6.4 of [LM].
\begin{lem}
Let $\pi$ be a permutation of a nonempty finite set $F$ and $\mathcal{C}=\{ E \subset F : \pi(E)=E \}$. Then 
$$\sum_{E \in \mathcal{C} \setminus \{ \varnothing \}} (-1)^{|E|+1} \, \mathrm{sgn} (\pi|_E)=1.$$
\end{lem}
With the help of this lemma, we obtain
\begin{equation}
p_m(\sigma_X) = \sum_{x \in P(m)} \sum_{E \in \mathcal{C}(m;\,x)\setminus \{ \varnothing \} } (-1)^{|E|+1} \text{sgn}(\sigma_A^m|_E).
\end{equation}
For $k=1,2,\cdots,|\mathcal{A}|$ let $X_k$ denote the topological Markov chain determined by the zero-one $\mathcal{A}_k \times \mathcal{A}_k$ matrix $B_k$:
$$X_k=\{ z\in\mathcal{A}_k^{\mathbb{Z}} : \forall i \in \mathbb{Z} \; F(z_i, z_{i+1}) \neq \varnothing \}.$$
We denote the shift map of $X_k$ by $\sigma_k$ and put
$$P(m; \, k)= \{ z \in X_k : \sigma_k^m(z)=z \}.$$
If $z \in P(m ;\, k)$, then $|F(z_i, z_{i+1})|=1$ for all $i$, because $\mathcal{L}_{\infty}: \textsf{X}_A \rightarrow X$ has no graph diamonds; 
hence $A_k(z_i, z_{i+1})=1$ or $-1$ for all $i$. 
In this case, we write
$$\text{sgn}(z)=A_k(z_0, z_1)A_k(z_1, z_2) \cdots A_k(z_{m-1}, z_m).$$
With this notation, we have
\begin{equation}
\sum_{z \in P(m;\, k)} \text{sgn}(z)=\text{tr}(A_k^m) \qquad(k=1, 2, \cdots, |\mathcal{A}|).
\end{equation}

Let $P_0$ denote the map $\textsf{X}_A \ni y \mapsto y_0 \in \mathcal{A}$. 
Suppose $x \in P(m)$.
Since $\mathcal{L}_{\infty}:\textsf{X}_A \rightarrow X$ has no graph diamonds and $x$ is periodic, 
the restriction of $P_0$ to $\mathcal{L}_{\infty}^{-1}(x)$ is one-to-one.
In particular, $|\mathcal{L}_{\infty}^{-1}(x)| \leq |\mathcal{A}|$.
If $E \in \mathcal{C}(m; \, x) \setminus \{ \varnothing\}$ and $k=|E|$, 
then $1 \leq k \leq |\mathcal{A}|,\, P_0(\sigma_A^i(E)) \in \mathcal{A}_k$ for all $i$ and
$$\big(P_0(\sigma_A^i(E))\big)_{i \in \mathbb{Z}} \in P(m; \, k).$$
We denote the point $(P_0(\sigma_A^i(E)))_{i \in \mathbb{Z}}$ by $\Phi(E)$.
With this notation, we have
\begin{equation}
\text{sgn}(\sigma_A^m|_E)=\text{sgn}(\Phi(E)), \qquad(E \in \mathcal{C}(m,\,x)\setminus \{ \varnothing \}).
\end{equation}
Finally, it is easy to check that the map
$$\Phi: \bigcup_{x \in P(m)} \mathcal{C}(m; \, x) \setminus \{ \varnothing \} \rightarrow \bigcup_{k=1}^{|\mathcal{A}|} P(m;\,k)$$
is a one-to-one correspondence.
Hence (2.4) follows from(2.5), (2.6)and(2.7).

In our proof of the other formulas in Theorem B, we need the following modification of Lemma 2.2.
\begin{lem}
Let $\pi$ be a permutation of a finite set $F$. Suppose $\mathcal{C} \subset 2^{F}$ and $G \subset F$ satisfy the following:
\newline
\textsc{(1)} $\pi(E)=E$ for every $E \in \mathcal{C}$.
\newline
\textsc{(2)} If $E_1, E_2 \in \mathcal{C}$, then $E_1 \cup E_2, E_1 \setminus E_2 \in \mathcal{C}$.
\newline
\textsc{(3)} $G \neq \varnothing$ and $\pi(G)=G$.
\newline
\textsc{(4)} If $E \subset G$ and $\pi(E)=E$, then $E \in \mathcal{C}$.
\newline
Then 
\begin{equation}
\sum_{E \in \mathcal{C} \setminus \{ \varnothing \}} (-1)^{|E|+1} \, \mathrm{sgn}(\pi|_E)=1. 
\end{equation}
\end{lem}
\begin{pf}
The assumptions imply that $\mathcal{C}$ has at least two distinct elements, namely $G$ and $\varnothing$. 
For convenience, we set $\text{sgn}(\pi|_{\varnothing})=1$. 
With this notation, (2.8) is equivalent to
$$\sum_{E \in \mathcal{C}} \, (-1)^{|E|}\, \text{sgn} (\pi|_E)=0.$$

If we put 
$$\mathcal{N}=\{ E \in \mathcal{C} : E \subset G\} \;\; \text{and} \;\; \mathcal{R} = \{ E\in \mathcal{C} : E \cap G = \varnothing \},$$
then (2) implies that
$$\mathcal{N} \times \mathcal{R} \ni (E_1, E_2) \mapsto E_1 \cup E_2 \in \mathcal{C}$$
is a one-to-one correspondence.
If $(E_1, E_2) \in \mathcal{N} \times \mathcal{R}$, then $E_1 \cap E_2 = \varnothing$,
and hence
$$(-1)^{|E_1 \cup E_2|}\, \text{sgn} (\pi|_{E_1 \cup E_2}) = (-1)^{|E_1|}\, \text{sgn}(\pi|_{E_1})\, (-1) ^{|E_2|}\, \text{sgn}(\pi|_{E_2}).$$

From (3) and Lemma 2.2, we have
$$\sum_{E \subset G, \, \pi(E)=E } (-1)^{|E|}\, \text{sgn}(\pi|_E)=0;$$
and from (1) and (4), we have
$$\{ E \subset G : \pi(E) =E \} = \mathcal{N}.$$
Therefore
\begin{eqnarray*}
\sum_{E\in \mathcal{C}} \, (-1)^{|E|} \, \text{sgn}(\pi|_E) 
&=& \sum_{(E_1, E_2) \in \mathcal{N} \times \mathcal{R}} (-1)^{|E_1 \cup E_2|} \, \text{sgn}(\pi|_{E_1 \cup E_2})\\
&=& \sum_{E_1 \in \mathcal{N}} \, (-1)^{|E_1|} \, \text{sgn}(\pi|_{E_1}) \sum_{E_2 \in \mathcal{R}} \, (-1)^{|E_2|} \, \text{sgn} (\pi|_{E_2})\\
&=&0. 
\end{eqnarray*}
\hfill$\Box$
\end{pf}
\begin{rmk}
The condition that $G \neq \varnothing$ is crucial. Otherwise we cannot apply Lemma 2.2.
\end{rmk}

Now, suppose that $N$ is a positive integer and $\delta \in \{0, 1\}$. We put
$$P(N, \delta)=\{ x \in X : \sigma_X^N(x) = \sigma_X^{\delta}\, \varphi(x)=x \}$$
so that $p_{N, \delta}(\sigma_x, \varphi)=|P(N, \delta)|$; 
and for $x \in P(N, \delta)$ we put
$$\mathcal{C}(N, \delta;\, x)=\{ E \subset \mathcal{L}_{\infty}^{-1}(x) : \sigma_A^N(E) = \sigma_A^{\delta} \, \varphi_{J,A}(E)=E \}.$$
We have $P(N, \delta)\subset P(N)$ and for $x\in P(N, \delta)$ we have $\mathcal{C}(N, \delta; \, x) \subset \mathcal{C}(N ;\, x)$.

\begin{lem}
We have
$$p_{N, \delta}(\sigma_X, \varphi)
=\sum_{x \in P(N, \delta)} \sum_{E\in \mathcal{C}(N, \delta ;\, x) \setminus \{ \varnothing\}} (-1)^{|E|+1} \, \mathrm{sgn}(\sigma_k^N|_E).$$
\end{lem}
\begin{pf}
Let $x \in P(N, \delta)$ and put
$$G(x)=\{ y \in \mathcal{L}_{\infty}^{-1}(x) : \exists n \in \mathbb{Z} \;\; \sigma_A^{\delta}\, \varphi_{J,A}(y) =\sigma_A^{nN}(y) \}.$$
It follows from (5) in Proposition 2.1 that $G(x) \neq \varnothing$.
It is then straightforward to check that the conditions in Lemma 2.3 are satisfied with 
$F=\mathcal{L}_{\infty}^{-1}(x), \, \pi=\sigma_A^N|_{\mathcal{L}_{\infty}^{-1}(x)}, \, \mathcal{C}=\mathcal{C}(N, \delta ;\, x)$ and
$G=G(x)$.
Hence we obtain
$$\sum_{E\in \mathcal{C}(N, \delta;\, x) \setminus \{\varnothing \}} (-1)^{|E|+1} \, \text{sgn} (\sigma_A^N|_E)=1 \qquad(x\in P(N, \delta)),$$
and the result follows. \hfill$\Box$
\end{pf}

Suppose $1\leq k \leq |\mathcal{A}|$ and let us consider the topological Markov chain $X_k$ again.
If $S \in \mathcal{A}_k$, we have $\tau(S) \in \mathcal{A}_k$.
We also have
$$B_k(S_1, S_2)=B_k(\tau(S_2), \tau(S_1)) \qquad(S_1, S_2 \in \mathcal{A}_k).$$
Hence the map $\varphi_k : X_k \rightarrow X_k$ defined by
$$\varphi_k(z)_i=\tau(z_{-i})$$
is a flip for $(X_k, \sigma_k)$.
We put
$$P(N, \delta;\, k) = \{ z \in X_k : \sigma_k^N(z)=\sigma_k^{\delta} \, \varphi_k(z)=z \}.$$

Again we have $P(N, \delta ; \, k) \subset P(N;\, k)$ for every $k$.
Moreover, if $x\in P(N, \delta), \, E \in \mathcal{C}(N, \delta ;\, x) \setminus \{ \varnothing \}$ and $k=|E|$,
then $1 \leq k \leq |\mathcal{A}|$ and $\Phi(E) \in P(N, \delta ; \, k)$;
and, as before, the map
$$\Phi : \bigcup_{x\in P(N, \delta)} \mathcal{C}(N, \delta ; \, x)\setminus \{ \varnothing \} \rightarrow \bigcup_{k=1}^{|\mathcal{A}|} P(N, \delta ;\, k)$$
is a one-to-one correspondence.

If $E \in \mathcal{C}(N, \delta ; \, x)\setminus \{ \varnothing \}$ for some $x\in P(N, \delta)$, 
then (2.7) implies that
$$\text{sgn}(\sigma_A^N|_E)=\text{sgn}(\Phi(E)).$$
By Lemma 2.4, we have
$$p_{N, \delta}(\sigma_X, \varphi) = \sum_{k=1}^{|\mathcal{A}|} \, (-1)^{k+1} \, \sum_{z \in P(N, \delta ; \, k )} \text{sgn}(z).$$
Therefore the following lemma completes our proof of Theorem B.
\begin{lem}
Let $1\leq k \leq |\mathcal{A}|$. Then we have
$$\sum_{z\in P(N, \delta ; \, k)} \mathrm{sgn}(z) = 
\begin{cases} 
\mathcal{S}\big[J_k^{\Delta} B_k^m J_k^{\Delta}\big] \qquad \qquad \qquad \; \; \text{if} \; N=2m \; \text{and} \; \delta=0, \\
\mathcal{S}\big[(J_k A_k)^{\Delta} B_k^{m-1} (A_k J_k)^{\Delta}\big] \quad \text{if} \; N=2m \; \text{and} \; \delta=1, \\
\mathcal{S}\big[J_k^{\Delta} B_k^m (A_k J_k)^{\Delta}\big] \qquad \quad \quad \; \text{if} \; N=2m+1 \; \text{and} \; \delta=0.  
\end{cases} $$
\end{lem}  
\begin{pf}
We prove the last one only.
The others are similarly proved.
Let $m$ be a non-negative integer.
We have to show that
$$\sum_{z\in P(2m+1, 0 ;\, k)} \text{sgn}(z) = S\big[J_k^{\Delta} B_k^m (A_k J_k)^{\Delta}\big].$$

For $S\in \mathcal{A}_k$ we write $\tau(S)=S^*$.
Thus $J_k(S_1, S_2) \neq 0$ if and only if $S_1 ^* =S_2$.

Let $S_0, S_1, \cdots,  S_m \in \mathcal{A}_k^{m+1}$.
We denote the product
$$J_k(S_0, S_0) B_k(S_0, S_1) \cdots B_k(S_{m-1}, S_m) A_k(S_m, S_m^*)J_k(S_m^*, S_m)$$
by $\Pi(S_0S_1 \cdots S_m)$.
If $m=0$,
$$\Pi(S_0)=J_k(S_0, S_0) A_k(S_0, S_0^*)J_k(S_0^*, S_0).$$
The expression $\mathbb{E}(S_0S_1 \cdots S_m)$ will denote the unique point $z\in \mathcal{A}_k^{\mathbb{Z}}$ such that
$$z_0 z_1 \cdots z_{2m} = S_0S_1 \cdots S_m S_m^* \cdots S_1^*$$
and
$$z_{i+2m+1}=z_i \qquad (i \in \mathbb{Z}).$$
If $m=0$, then $\mathbb{E}(S_0)$ is the point in $\mathcal{A}_k^{\mathbb{Z}}$ such that $z_i=S_0$ for all $i$.

We have
$$S\big[J_k^{\Delta} B_k^m (A_k J_k)^{\Delta}\big] = \sum_{S_0S_1 \cdots S_m \in \mathcal{I}} \Pi(S_0S_1 \cdots S_m),$$
where
$$\mathcal{I} = \{ S_0S_1 \cdots S_m \in \mathcal{A}_k^{m+1} : \Pi(S_0S_1 \cdots S_m) \neq 0 \},$$
and it is straightforward to check that the map
$$\mathcal{I} \ni S_0S_1 \cdots S_m \mapsto \mathbb{E}(S_0S_1 \cdots S_m) \in P(2m+1 , 0 ;\, k)$$
is well-defined, one-to-one and onto.
It remains only to show that
$$\Pi(S_0S_1 \cdots S_m) = \text{sgn}(\mathbb{E}(S_0S_1 \cdots S_m))$$
holds for every $S_0S_1 \cdots S_m \in \mathcal{I}$.

Let $S_0S_1 \cdots S_m \in \mathcal{I}$.
We have to show that
\begin{eqnarray}
&&\Pi(S_0S_1 \cdots S_m) \\
&&= A_k (S_0, S_1) \cdots A_k(S_{m-1}, S_m)A_k(S_m, S_m^*)\cdots A_k(S_1^*, S_0).\nonumber
\end{eqnarray}
If we put $z=\mathbb{E}(S_0S_1 \cdots S_m)$, then $|A_k(z_i, z_{i+1})|=1$ for all $i$, and hence
\begin{eqnarray*}
&&\Pi(S_0S_1 \cdots S_m) \\
&&= J_k(S_0, S_0)\big(A_k (S_0, S_1) \cdots A_k(S_{m-1}, S_m)\big)^2 A_k(S_m, S_m^*)J_k(S_m^*, S_m).
\end{eqnarray*}
We also have $S_0^* =S_0$, because $J_k(S_0, S_0) \neq 0$.

It is easy to show that
$$J_k(S, \mathcal{T})=J_k(\mathcal{T}, S) \;\; \text{and} \;\; J_k(S, S^*)A_k(S, \mathcal{T})=A_k(\mathcal{T}^*, S^*)J_k(\mathcal{T}, \mathcal{T}^*)$$
holds for all $S, \mathcal{T} \in \mathcal{A}_k$.
Thus
\begin{eqnarray*}
&& J_k(S_0, S_0)A_k(S_0, S_1) \cdots A_k(S_{m-1}, S_m) J_k(S_m^*, S_m) \\
&=& A_k(S_1^*, S_0^*) \cdots A_k(S_m^*, S_{m-1}^*)J_k(S_m, S_m^*)J_k(S_m^*, S_m)\\
&=& A_k(S_m^*, S_{m-1}^*)\cdots A_k(S_1^*, S_0).
\end{eqnarray*}
Hence (2.9) follows. \hfill $\Box$
\end{pf}
\section{Krieger's Joint State Chain(Proof of Proposition 2.1)}
We start with some notation.
Suppose $\mathcal{A}$ is a finite alphabet, $j, k \in \mathbb{Z}$ and $j \leq k$.
If $x \in \mathcal{A}^{\mathbb{Z}}$, the expression $x_{[j,\,k]}$ will denote the restriction of $x$ to the set
$[j,\, k] = \{ i\in \mathbb{Z} : j \leq i \leq k \}$.
Thus $x_{[j,\, k]} \in \mathcal{A}^{[j,\,k]}$.
If $E \subset \mathcal{A}^{\mathbb{Z}}$, we also write $E_{[j, \,k]} = \{ x_{[j, \, k]} : x \in E \}$.
We will use similar notation for other subsets of $\mathbb{Z}$ as well.

The following lemma implies that every flip for a shift dynamical system is conjugate to a one-block one.
\begin{lem}
Suppose $X$ is a shift space and $\varphi$ is a flip for $(X,\, \sigma_X)$.
Then there are a finite set $\mathcal{A}$, a map $\Phi : \mathcal{A} \rightarrow \mathcal{B}_1(X)$, 
a shift space $Y$ over $\mathcal{A}$ and a one-block flip $\psi$ for $(Y,\, \sigma_Y)$ such that $\Phi_{\infty}$ is a conjugacy from 
$(Y,\, \sigma_Y, \, \psi)$ to $(X, \, \sigma_x, \, \varphi)$.
\end{lem}
\begin{pf}
Since $\varphi$ is (uniformly) continuous, there is a non-negative integer $N$ such that
$$x, x' \in X \;\; \text{and} \;\; x_{[-N,\,N]} = x'_{[-N, \, N]}\; \Rightarrow\; \varphi(x)_0 = \varphi(x')_0 .$$
Put
$$\mathcal{A}=\{ (a, b) :  \exists x\in X \;\; \text{s.t.}\;\; a=x_0 \;\;\text{and}\;\; b=\varphi(x)_0 \}$$
and define the maps $\Phi : \mathcal{A} \rightarrow \mathcal{B}_1(X)$ and $\tau: \mathcal{A} \rightarrow \mathcal{A}$ by
$$\Phi(a, b)=a \;\; \text{and} \;\; \tau(a, b)=(b, a).$$
It is clear that $\mathcal{A}$ is a finite set, $\Phi$ and $\tau$ are well-defined, and that $\tau^2=\text{id}_{\mathcal{A}}$.

Define $\theta:X\rightarrow \mathcal{A}^{\mathbb{Z}}$ by
$$\theta(x)_i = \big( x_i \, , \; \varphi(x)_{-i} \big).$$
Since $\varphi$ is (uniformly) continuous, there is a non-negative integer $N$ such that
$$x, x' \in X \;\; \text{and} \;\; x_{[-N,\,N]} = x'_{[-N, \, N]}\; \Rightarrow\; \varphi(x)_0 = \varphi(x')_0 .$$
So, $\theta$ is a sliding block code with memory $N$ and anticipation $N$,
and $\theta(X)$ is a shift space over $\mathcal{A}$.

We set $Y=\theta(X)$ and define $\psi : Y \rightarrow Y$ by
$$\psi(y)_i = \tau (y_{-i}).$$
It is then straightforward to check that $\psi$ is well-defined, is a flip for $(Y, \, \sigma_Y)$
and that $\Phi_{\infty}$ is a conjugacy from $(Y, \, \sigma_Y, \, \psi)$ to $(X, \, \sigma_X, \, \varphi)$.
\hfill $\Box$
\end{pf}

Now, we prove Proposition 2.1.
Let $X$ be a sofic shift and $\varphi$ a flip for $(X, \, \sigma_X)$.
With the aid of Lemma 3.1, we may assume that $\varphi$ is a one-block flip.
We put $\mathcal{B}_1(X)=\mathcal{B}$ and denote the symbol map of $\varphi$ by $\tau$.
Thus $\tau : \mathcal{B} \rightarrow \mathcal{B}, \; \tau^2=\text{id}_{\mathcal{B}}$ and we have
$$\varphi(x)_i = \tau(x_{-i}) \qquad(x \in X , \; i \in \mathbb{Z}).$$
We prove Proposition 2.1 by showing that Krieger's joint state chain of the sofic system $(X, \sigma_X)$ has the desired properties.
In the case when $X$ is irreducible(topologically transitive) are further show that Krieger's joint finitary state chain also has the desired properties.

We denote the set of all left-infinite sequences, right-infinite sequences and blocks in $X$ by $L$, $R$ and $B$, respectively:
$$L=\bigcup_{j\in \mathbb{Z}}X_{(-\infty, \, j]}, \; R=\bigcup_{j\in \mathbb{Z}} X_{[j,\, \infty)} \;\; \text{and} \;\; B=\bigcup_{-\infty < j \leq k < \infty} X_{[j, \, k]}.$$
The rules ``$\sigma(x)_i = x_{i+1}$'' and ``$\varphi(x)_i=\tau(x_{-i})$'' will be applied to those sequences as well.
For instance, if $\lambda \in X_{(-\infty,\, j]}$, then $\sigma(\lambda)\in X_{(-\infty,\,j-1]}, \; \varphi(\lambda)\in X_{[-j, \, \infty)}$,
$$\sigma(\lambda)_i=\lambda_{i+1}\qquad (i \leq j-1),$$
and
$$\varphi(\lambda)_i=\tau(\lambda_{-i}) \qquad (i \geq -j).$$
For notational simplicity, we write $\varphi(x)=x^*$ whenever $x\in L\cup R \cup B$.
Thus we have $(x^*)^*=x$ for all $x$,
$$L=\{ \rho^* :\rho \in R  \}\,,\, R=\{ \lambda^* : \lambda \in L \}\;\; \text{and} \;\; B=\{ w^* : w \in B \}.$$
We also write $a^*=\tau(a)$ for $a \in \mathcal{B}$.

If $j, k \in \mathbb{Z}, \, \lambda \in X_{(-\infty,\, j]}$ and $\rho \in X_{[k, \infty)}$, we denote the bi-infinite sequence $x \in \mathcal{B}^{\mathbb{Z}}$ 
such that $x_{(-\infty,\,-1]}=\sigma^{j+1}(\lambda)$ and $x_{[0, \, \infty)} = \sigma^k(\rho)$ by $\lambda.\rho\;$:
$$(\lambda.\rho)_i = \begin{cases} \lambda_{i+j+1} \qquad(i\leq -1), \\ \rho_{i+k} \qquad \quad (i \geq 0).\end{cases}$$
We will use similar notation for other sequences in $L \cup R \cup B$ as well.

If $\lambda \in L$, its future is the set $\mathbb{F}(\lambda) \subset X_{[0, \, \infty)}$ defined by
$$\mathbb{F}(\lambda) = \{ \rho \in X_{[0, \infty)} : \lambda . \rho \in X \}.$$
The pasts are similarly defined:
$$\mathbb{P}(\rho) = \{ \lambda \in X_{(-\infty, \, 0]} : \lambda . \rho \in X \} \qquad(\rho \in R).$$
Is is well known that the sets
$$\varXi_+ = \{ \mathbb{F}(\lambda) : \lambda \in L \} \;\; \text{and} \;\; \varXi_- = \{ \mathbb{P}(\rho) : \rho \in R \} $$
are finite [K].
A triple $(F, a, P)\in \varXi_+ \times \mathcal{B} \times \varXi_-$ is called a joint state of $(X, \, \sigma_X)$
if $a\in P_0(F) \cap P_0(P)$, where
$$P_0(F)=\{ \rho_0 : \rho \in F \} \; \text{and} \; P_0(P)=\{ \lambda_0 : \lambda \in P \}.$$
We denote the set of all joint states of $(X, \, \sigma_X)$ by $\Omega$ and define
$\mathcal{L} : \Omega \rightarrow \mathcal{B}$ by 
$$\mathcal{L}(F, a, P)=a.$$

If $F\in \varXi_+$ and $a \in P_0(F)$, then the set
$$F(a) = \{ \rho \in X_{[0, \, \infty)} : \sigma^{-1}(a.\rho) \in F \}$$
is an elements of $\varXi_+$; and $P \in \varXi_-$ and $a \in P_0(P)$, then the set
$$P(a) = \{ \lambda \in X_{(-\infty,\, 0]} : \lambda.a \in P \}$$
is an element of $\varXi_-$.
We define the zero-one $\Omega \times \Omega$ matrix $A$ by
$$A((F_1, a_1, P_1), (F_2, a_2, P_2)) = 1 \;\; \Leftrightarrow \;\; F_2 = F_1(a_1)\;\; \text{and}\;\; P_1=P_2(a_2).$$
The topological Markov chain $\textsf{X}_A$ is called Krieger's joint state chain of $(X, \, \sigma_X).$
It is well known that $\mathcal{L}_{\infty} : \textsf{X}_A \rightarrow X$ is a factoring [K]
and it is easy to check that $\mathcal{L}_{\infty}$ has no graph diamonds.

If $F \in \varXi_+$ and $P\in \varXi_-$, then the sets
$$F^* = \{ \rho^* : \rho \in F \}\;\;  \text{and} \;\; P^* = \{ \lambda^* : \lambda \in P \}$$
are elements of $\varXi_-$ and $\varXi_+$, respectively.
For $(F,\, a,\, P)\in \Omega$ we write
$$(F,\, a,\, P)^* = (P^*, \,a^*,\, F^*).$$
If $\mathbf{a} \in \Omega$, we have $\mathbf{a}^* \in \Omega$ and $(\mathbf{a}^*)^*=\mathbf{a}$.
We also have
$$A(\mathbf{a}, \mathbf{b}) = A(\mathbf{b}^*, \, \mathbf{a}^*) \qquad(\mathbf{a}, \mathbf{b} \in \Omega).$$
Consequently, if we put $\mathcal{A} = \Omega$ and define the zero-one $\mathcal{A} \times \mathcal{A}$ matrix $J$ by
$$J(\mathbf{a}, \mathbf{b}) = 1 \;\;  \Leftrightarrow \;\; \mathbf{a}^* = \mathbf{b},$$
then the system $(\mathcal{A}, \mathcal{L}, A, J)$ satisfies the condition (1) through (4) in Proposition 2.1.

To prove (5), suppose $x\in X$.
Since $\mathcal{L}_{\infty}: \textsf{X}_A \rightarrow X$ is a factoring, there is a $y \in \textsf{X}_A$ such that
$\mathcal{L}_{\infty}(y) = x$.
We may write
$$y = (F_i, \, x_i,\,  P_i)_{i\in \mathbb{Z}}.$$
We then have 
$$\varphi_{J,A}(y)_i = (P_{-i}^*,\, x_{-i}^*,\, F_{-i}^*) \qquad(i \in \mathbb{Z})$$
and
$$\sigma_A \,\varphi_{J,A}(y) = (P_{-(i+1)}^*,\, x_{-(i+1)}^*, \, F_{-(i+1)}^*) \qquad (i \in \mathbb{Z}).$$
If $\varphi(x)=x$, then $x_i = (x_{-i})^*$ for all $i$ and hence the point
$$y'=(F_i, \,x_i,\, F_{-i}^*)_{i \in \mathbb{Z}}\in \textsf{X}_A$$
satisfies $\mathcal{L}_{\infty}(y')=x$ and $\varphi_{J,A}(y')=y'$; and if $\sigma_X \, \varphi(x)=x$, 
then $x_i=(x_{-(i+1)})^*$ for all $i$ and hence the point
$$y''=(F_i, \,x_i,\, F_{-(i+1)}^*) \in \textsf{X}_A$$
satisfies $\mathcal{L}_{\infty}(y'')=x$ and $\sigma_A \, \varphi_{J,A}(y'')=y''$.
This completes the proof of Proposition 2.1.

A block $w \in B$ is said to be finitary(intrinsically synchronizing) if it has following property:
\newline
If $\lambda \in L, \, \rho \in R, \, \lambda.w \in L$ and $w.\rho \in R$, then $(\lambda.w).\rho \in X$.
\newline  
The set of all finitary blocks in $X$ will be denoted by $B^0$. 
It is obvious that if $w \in B^0$, then the sets
$$\mathbb{F}(w)=\{ \rho \in X_{[0, \, \infty)} : w.\rho \in R \} \;\; \text{and} \;\; 
\mathbb{P}(w) = \{ \lambda \in X_{(-\infty, \, 0]} : \lambda . w \in L \}$$
are elements of $\varXi_+$ and $\varXi_-$, respectively.
It is also obvious that $w \in B^0$ if and only if $w^* \in B^0$.

We put
$$\varXi_+^0=\{ \mathbb{F}(w) : w \in B^0 \}, \;\; \varXi_-^0 = \{ \mathbb{P}(w) : w \in B^0 \}$$
and
$$\Omega^0 = \{ (F, a, P) \in \Omega : F \in \varXi_+^0 \;\; \text{and} \;\; P \in \varXi_-^0 \}.$$
We denote the resrictions of $\mathcal{L}, \; A$ and $J$ to $\Omega^0$ by $\mathcal{L}^0, \; A^0$ and $J^0$, respectively.
The topological Markov chain $\textsf{X}_{A^0}$ is called Krieger's joint finitary state chain of $(X, \, \sigma_X)$.

If $X$ is irreducible, it is well known that $(\mathcal{L}^0)_{\infty} : \textsf{X}_{A^0} \rightarrow X$ is a factoring [K].
Now, it is clear that if $X$ is irreducible, then the system $(\Omega^0, \mathcal{L}^0, A^0, J^0)$ also satisfies the conditions in Proposition 2.1.
\section{An Example and Some Remarks}
\begin{rmk}
Let $X$ be a sofic shift and $\varphi$ a one-block flip for $(X, \sigma_{X})$.
Suppose $(\mathcal{A}, \mathcal{L}, A, J)$ satisfy the conditons (1) through (4) in Proposition 2.1.
Then the conditon (5) is equivalent to the followings:
\newline
(i) for any $a \in \mathcal{B}_1(X)$ and $\rho \in F(a)$ with $a^*=a$ and $\rho^{*}.a\rho \in X$ 
there exist $b\in {\mathcal{L}}^{-1}(a)$ and $\lambda \in {\mathcal{L}}^{-1}(\rho)$ such that $b^*=b$ and $\sigma^{-1}(b.\lambda) \in (\textsf{X}_A)_{[0, \infty]}$, and
\newline
(ii) for any $a \in \mathcal{B}_1(X)$ and $\rho \in F(a)$ with $a^*a \in \mathcal{B}_2(X)$ and $\rho^{*}a^*.a\rho \in X$ 
there exist $b\in {\mathcal{L}}^{-1}(a)$ and $\lambda \in {\mathcal{L}}^{-1}(\rho)$ such that $b^*b \in \mathcal{B}_2(\textsf{X}_A)$ and $\sigma^{-1}(b.\lambda) \in (\textsf{X}_A)_{[0, \infty]}$.
\end{rmk}
The joint finitary state chain of a irreducible sofic shift may not be irreducible.
If we take
$$\Omega' = \{ (\mathbb{F}(w_1), \; a, \; \mathbb{P}(w_2) )\in \Omega^0 : w_1 a w_2 \in B \}, $$
then $(\Omega', \mathcal{L}', A', J')$ is a irreducible component which satisfies the conditions (1) through (4) in Proposition 2.1,
where $\mathcal{L}, \; A'$ and $J'$ are restrictions of $\mathcal{L}, \; A$ and $J$ to $\Omega'$.
However, it may not satisfy the condition (5).

For example, the even shift $(X, \sigma)$ with the flip map $\varphi$ induced by the symbol map $\tau=\text{id}_{\{0,1\}}$ has a fixed point
$0^{\infty} \in P(2m, 0)$ for any positive integer $m$, but there is not $y \in (\textsf{X}_{A'}, \sigma_{A'})$ such that 
${\sigma_{A'}}^{2m}(y)=\varphi_{J'A'}(y)=y$. In this case, $(\Omega', \mathcal{L}', A', J')$ dose not satisfy the condition (i) of the above Remark.

\begin{exa}
Let $(X, \sigma)$ be the even shift, and $\varphi:X \rightarrow X$ be the flip induced by the symbol map $\tau=\mathrm{id}_{\{ 0, 1 \}}$.
Then the zeta fuction is given by
$$\zeta_{\sigma, \varphi}(t)=\sqrt{\frac{1+t^2}{1-t^2-t^4}}\, \exp\left(\frac{2t+2t^2-t^3-t^4-2t^5-t^6}{(1-t^2)(1-t^2-t^4)}\right).$$
\end{exa}

\begin{figure}[!htb]
\centering
\includegraphics[scale=0.5]{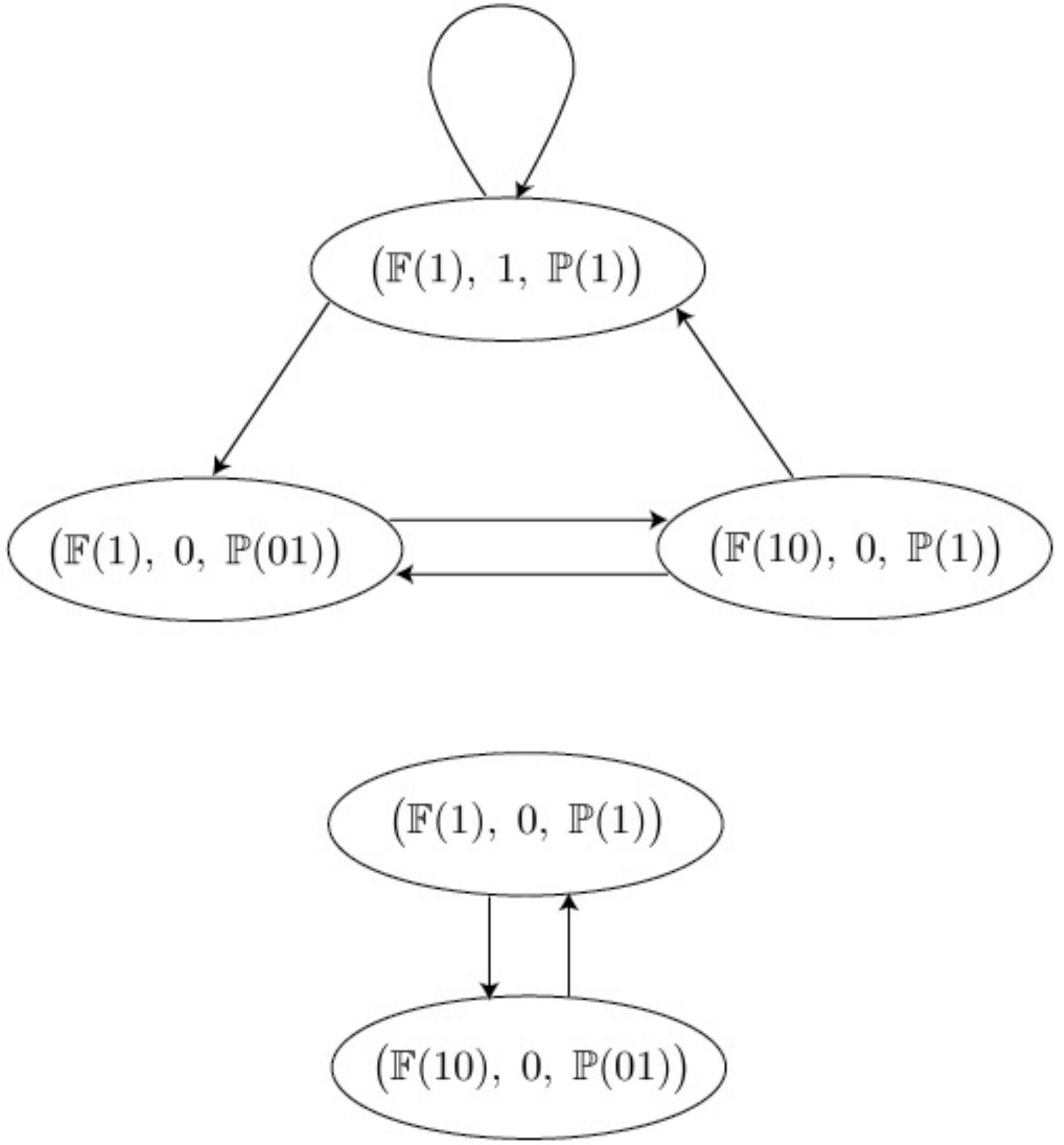}
\\
Figure. The joint finitary state chain of the even shift
\end{figure}
 
\begin{center}
\textbf{References}
\end{center}
\vspace{0.3cm}

\noindent
[B] R. Bowen, \textit{On Axiom A Diffeomorphisms}, AMS-CBMS Reg. Conf. 35, Providence, 1978
\vspace{0.3cm}
\newline
[K] W. Krieger, \textit{On Sofic Systems I}, Israel Journal of Mathematics, 48(1984), 305-332.
\vspace{0.3cm}
\newline
[KLP] Y.-O. Kim, J. Lee and K. K. Park, \textit{A Zeta Function for Flip Systems}, Pacific Journal of Mathematics, 209(2003), 289-301.
\vspace{0.3cm}
\newline
[L] D. Lind, \textit{A zeta function for $\mathbb{Z}^d$-actions}, London Math. Soc. Lecture Note Ser., 228
(M. Pollicott and K.Schmidt, ed.), 1996, 433-450.
\vspace{0.3cm}
\newline
[LM] D. Lind and B. Marcus, \textit{Symbolic Dynamics and Coding}, Cambridge University Press, 1995.
\vspace{0.3cm}
\newline
[M] A. Manning, \textit{Axiom A diffeomorphisms have rational zeta functions}, Bull. London Math. Soc. 3(1971), 215-220

\end{document}